\documentclass[12pt]{amsart}

\usepackage{latexsym, amssymb, amscd, amsfonts}

\usepackage[colorlinks=true,urlcolor=black,citecolor=black,linkcolor=black,backref=section,%
pdftitle={Intersection multiplicity one for classical groups},%
pdfauthor={Ivan Dimitrov and Mike Roth},%
pdfsubject={Geometry and representation theory},%
pdfkeywords={Algebraic geometry, Cohomology of Homogeneous Spaces, Representation theory}]{hyperref}
\usepackage{ifthen}
\usepackage{longtable}

\DeclareMathAlphabet{\mathpzc}{OT1}{pzc}{m}{it}

\newboolean{showpicture}
\setboolean{showpicture}{true}

\newboolean{bourbaki}
\setboolean{bourbaki}{true}

\newboolean{draft}
\setboolean{draft}{true}

\newboolean{shproofs}
\setboolean{shproofs}{false}

\newcommand{\np}{\medskip\noindent}

\newcommand{\point}{\vspace{3mm}\par\refstepcounter{subsection}\noindent{\bf \thesubsection.} }
\newcommand{\tpoint}[1]{\vspace{3mm}\par\refstepcounter{subsubsection}\noindent{\em #1 {\rm(}{\em \thesubsubsection}{\rm)} ---} }

\newcommand{\bpoint}[1]{\vspace{3mm}\par\refstepcounter{subsection}\noindent{\bf \thesubsection.} {\bf #1.} }

\renewenvironment{equation}{\medskip\noindent\refstepcounter{subsubsection}\makebox[0pt][l]{({\bf\thesubsubsection})}\begin{minipage}[b]{\textwidth}$$}{$$\end{minipage}\medskip\noindent}

\renewcommand{\labelenumi}{{({\em \alph{enumi}})}}

\ifthenelse{\boolean{draft}}
{
}
{
}

\newcommand{\bpf}{\noindent {\em Proof.  }}
\newcommand{\epf}{\qed \vspace{+10pt}}

\newcounter{enumcount}
\newcommand{\PauseEnumerate}{\end{enumerate}\setcounter{enumcount}{\arabic{enumi}}}
\newcommand{\ResumeEnumerate}{\begin{enumerate}\setcounter{enumi}{\theenumcount}}

\newcommand{\AlphaList}{\renewcommand{\labelenumi}{{({\em\alph{enumi}})}}}
\newcommand{\RomanList}{\renewcommand{\labelenumi}{{({\em\roman{enumi}})}}}

\renewcommand{\geq}{\geqslant}
\renewcommand{\leq}{\leqslant}

\newcommand{\st}{\,\,|\,\,}                            
\newcommand{\bst}{\,\,\,\rule[-0.3cm]{0.01cm}{0.8cm}\,\,\,}    

\renewcommand{\H}{\mathrm{H}}                          
\newcommand{\A}{\mathrm{A}}                            
\newcommand{\B}{\mathrm{B}}                            
\newcommand{\C}{\mathrm{C}}                            
\newcommand{\G}{\mathrm{G}}                            
\renewcommand{\P}{\mathrm{P}}                          
\newcommand{\X}{\mathrm{X}}                            
\newcommand{\M}{\mathrm{M}}                            
\newcommand{\N}{\mathrm{N}}                            
\newcommand{\U}{\mathrm{U}}                            
\newcommand{\Q}{\mathrm{Q}}                            
\newcommand{\W}{\mathcal{W}}                           
\newcommand{\WW}{\mathrm{W}}                           
\newcommand{\delpos}{\Delta^{\!+}}                     
\newcommand{\invset}{\Phi}                             
\newcommand{\wo}{w_{0}}                                
\renewcommand{\c}{\operatorname{c}}                    
\newcommand{\T}{\mathrm{T}}                            
\newcommand{\D}{\mathrm{D}}
\renewcommand{\SS}{\mathrm{S}}                         
\newcommand{\BKprod}{\odot_{0}}
\renewcommand{\emptyset}{\varnothing}
\newcommand{\sQ}{\mathcal{Q}}                          
\newcommand{\op}{\mbox{\scriptsize op}}                 

\newcommand{\gu}{\mathfrak{u}}
\newcommand{\gt}{\mathfrak{t}}

\newcommand{\taubar}{\underline{\tau}}

\newcommand{\ep}{\epsilon}                             

\newcommand{\codim}{\operatorname{codim}}

\newcommand{\Inv}{\Phi}

\newcommand{\rank}{\operatorname{rank}}                

\newcommand{\I}{\mathrm{I}}

\ifthenelse{\boolean{bourbaki}}
{
\newcommand{\CC}{\mathbf{C}} 
\newcommand{\PP}{\mathbf{P}} 
\newcommand{\QQ}{\mathbf{Q}} 
\newcommand{\ZZ}{\mathbf{Z}} 
}
{
\newcommand{\CC}{\mathbb{C}} 
\newcommand{\PP}{\mathbb{P}} 
\newcommand{\QQ}{\mathbb{Q}} 
\newcommand{\ZZ}{\mathbb{Z}} 
}

\ifthenelse{\boolean{draft}}
{

\newcommand{\remind}[1]{{\bf[#1]}}
\newcommand{\lremind}[1]{{\bf[label:  #1]}}
\newcommand{\bremind}[1]{{\bf[label:  #1]}}

\newcommand{\comment}[1]{{\bf [#1]}}
}
{

\newcommand{\remind}[1]{{}}
\newcommand{\lremind}[1]{{}}
\newcommand{\bremind}[1]{{}}

\newcommand{\comment}[1]{{}}
}

\newcommand{\hiddenproof}[1]{
\ifthenelse{\boolean{shproofs}}{
\medskip
\begin{centering}
\begin{minipage}{0.9\textwidth}
\hrule
\vspace{-0.75\baselineskip}
\small
#1
\vspace{0.25\baselineskip}
\hrule
\end{minipage}\\
\end{centering}
\medskip
}
{
}
}

\begin{document}
\pagestyle{plain} \title{{ \large{intersection Multiplicity One for Classical Groups}}
}
\author{Ivan Dimitrov}
\address{Department of Mathematics and Statistics, Queen's University, Kingston,
Ontario,  K7L 3N6, Canada} 
\email{dimitrov@queensu.ca} 
\thanks{Research of I.\ Dimitrov and M.\ Roth was partially supported by NSERC grants}
\author{Mike Roth}
\email{mike.roth@queensu.ca}


\begin{abstract}
In this paper we show that when $\G$ is a classical semi-simple algebraic group, $\B\subset\G$ a Borel
subgroup, and $\X=\G/\B$, then the structure coefficients of the Belkale-Kumar product $\BKprod$ on $\H^{*}(\X,\ZZ)$
are all either $0$ or $1$. 

\np
\end{abstract}

\maketitle

\vspace{0.5cm}

\section{Introduction}

\point Let $\G$ be a semi-simple algebraic group over an algebraically closed field of characteristic zero, 
$\B\subset \G$ a Borel subgroup, and set $\X=\G/\B$. 

\np
For any element $w$ of the Weyl group $\W$ of $\G$ the {\em Schubert variety} $\X_{w}$ is defined by 
$$\X_{w} := \overline{\B w\B/\B}\subseteq \G/\B=\X.$$
Recall that the classes of the Schubert cycles $\{[\X_{w}]\}_{w\in \W}$ give 
a basis for the cohomology ring $\H^{*}(\X,\ZZ)$ of $\X$.
Each $[\X_{w}]$ is a cycle of complex dimension $\ell(w)$, where $\ell(w)$ is the length of $w$.
The dual Schubert cycles $\{[\Omega_{w}]\}_{w\in \W}$, given by $\Omega_{w}:=\X_{\wo w}$, where $\wo\in \W$ is the
longest element, also form a basis. 
Each $[\Omega_{w}]$ is a cycle of complex codimension $\ell(w)$.

\np
For any $w_1,w_2,w \in \W$ we define the {\em structure constant} $c_{w_1,w_2}^{w}$
to be the coefficient of $[\Omega_{w}]$ when expressing the product $[\Omega_{w_1}]\cdot[\Omega_{w_2}]$ as a sum
of basis elements, so that
$$[\Omega_{w_1}]\cdot[\Omega_{w_2}] = \sum_{w\in \W} c_{w_1,w_2}^{w} [\Omega_{w}].$$
In \cite{BK} Belkale and Kumar define a new product $\BKprod$ on $\H^{*}(\X,\ZZ)$.
(More generally
\cite{BK} defines a new product on $\H^{*}(\G/\P,\ZZ)$, where $\P$ is any parabolic, however this paper is
only concerned with the case $\P=\B$.) 
Let $d_{w_1,w_2}^{w}$ be the structure coeffients of the 
Belkale-Kumar product, so that as above
$$[\Omega_{w_1}]\BKprod[\Omega_{w_2}] = \sum_{w\in \W} d_{w_1,w_2}^{w} [\Omega_{w}].$$
The Belkale-Kumar constants $d_{w_1,w_2}^{w}$ are equal to the usual constants $c_{w_1,w_2}^{w}$ if the triple 
$(w_1,w_2,w)$ is {\em Levi-movable} \cite[Definition 4]{BK}, and zero otherwise.   
Specifically, let $\Delta^{+}$ denote the set of positive roots
of $\G$, and $\Delta^{-}=-\Delta^{+}$ the negative roots.  Following Kostant \cite[Definition 5.10]{Ko}, 
for each $w\in \W$ we define the {\em inversion set} $\Inv_{w}:=w^{-1}\Delta^{-}\cap \Delta^{+}$.  
Belkale and Kumar \cite[Theorem 43+Corollary 44]{BK} prove that

\begin{equation}\label{eq:structure-coeffs}
d_{w_1,w_2}^{w} =\left\{
\begin{array}{cl}
c_{w_1,w_2}^{w} & \mbox{if $\Inv_{w}=\Inv_{w_1}\sqcup\Inv_{w_2}$} \\
0 & \mbox{otherwise,}\rule{0cm}{0.5cm} \\
\end{array}
\right.
\end{equation}
where $\sqcup$ denotes disjoint union.
The following two statements are therefore equivalent~:

\RomanList
\begin{enumerate}
\item The structure constants of the Belkale-Kumar product $\BKprod$ on $\H^{*}(\X,\ZZ)$ are all either $0$ or $1$.
\item $c_{w_1,w_2}^{w}=1$ whenever $\Inv_{w}=\Inv_{w_1}\sqcup \Inv_{w_2}$.
\PauseEnumerate
It is useful to write ({\em ii}) in a more symmetric form.
Since $\wo\Delta^{+}=\Delta^{-}$, 
it follows easily that $\Inv_{\wo w} = \Delta^{+}\setminus \Inv_{w}$, so that the condition 
$\Inv_{w}=\Inv_{w_1}\sqcup \Inv_{w_2}$ is equivalent to $\Delta^{+}=\Inv_{w_1}\sqcup \Inv_{w_2}\sqcup\Inv_{\wo w}$.
Furthermore, since the class $[\Omega_{\wo w}]$ is dual to $[\Omega_{w}]$ we have
$c_{w_1,w_2}^{w} = [\Omega_{w_1}]\cdot[\Omega_{w_2}]\cdot[\Omega_{\wo w}]$.  Setting $w_3=\wo w$ we can
therefore rephrase ({\em ii}) as 

\ResumeEnumerate
\item 
$[\Omega_{w_1}]\cdot[\Omega_{w_2}]\cdot[\Omega_{w_3}]=1$ whenever $\Delta^{+}=\Inv_{w_1}\sqcup\Inv_{w_2}\sqcup \Inv_{w_3}$.
\PauseEnumerate

\np
This in turn is equivalent to the following similar statement with an arbitrary number of elements of $\W$~:

\ResumeEnumerate
\item $\bigcap_{i=1}^{k}[\Omega_{w_i}] = 1$ whenever 
$w_1,\ldots, w_k\in\W$ satisfy 
\PauseEnumerate

\begin{equation}\label{eqn:lining-up}\label{eqn:sym-liningup}
\Delta^{+}=\bigsqcup_{i=1}^{k} \Inv_{w_i}.
\end{equation}

\np
It is clear that ({\em iv}) implies ({\em iii}).  The proof that ({\em iii}) implies ({\em iv}) requires
a slightly longer argument, and we defer it to the Appendix in \S\ref{sec:Appendix}.

\np
The main result of this paper is that these equivalent conditions hold for any classical group $\G$ and for the 
exceptional group $\G_2$, and hence for any semisimple group whose factors are of classical type or isomorphic 
to $\G_2$.
This result, proven in the form of ({\em iv})\footnote{Although, by the equivalences above,
it would suffice to prove only ({\em iii}), we have chosen to prove statement ({\em iv}) for arbitrary $k$ since 
it seems useful to record the more general versions of some of the combinatorial statements used 
in the proof.} appears as Theorem \ref{thm:mult1}.

\np
In the rest of this introduction we indicate a few other statements equivalent to the ones above.

\bpoint{Other equivalent statements}

\RomanList
\noindent
By \cite[Corollary 44]{BK} ({\em ii}) is equivalent to 
\ResumeEnumerate
\item 
$\prod_{\alpha\in \Inv_{w^{-1}}} \langle \rho,\alpha\rangle =
\left(\prod_{\alpha\in \Inv_{w_1^{-1}}} \langle \rho,\alpha\rangle \right)
\left(\prod_{\alpha\in \Inv_{w_2^{-1}}} \langle \rho,\alpha\rangle \right)
$
whenever $\Inv_{w}=\Inv_{w_1}\sqcup \Inv_{w_2}$.
\PauseEnumerate
Here $\rho$ is one-half the sum of the positive roots and $\langle\cdot,\cdot\rangle$ the Killing form. 

\np
In \cite[Theorem 43]{BK} Belkale and Kumar give an isomorphism of graded rings :
$$\phi\colon \left(\H^{*}(\X,\CC),\BKprod\right) \cong \left[\H^{*}(\gu^{+})\otimes \H^{*}(\gu^{-})\right]^{\gt},$$
where $\H^{*}(\gu^{\pm})$ denotes Lie algebra cohomology of the nilpotent algebras $\gu^{\pm}$, 
and $\gt$ the subalgebra corresponding to the maximal torus.  Under this isomorphism 
$$\phi\left([\Omega_{w}]\right) = 
\left(-1\right)^{\frac{p(p-1)}{2}}\left(\frac{i}{2\pi}\right)^{p}\left(\prod_{\alpha\in \Inv_{w^{-1}}}
\langle\rho,\alpha\rangle\right)\xi^{w}$$
where $p=\ell(w)$, and where (roughly) $\xi^{w}=\left(\wedge_{i=1}^{p} y_{\beta_i}\right)\otimes
\left(\wedge_{i=1}^{p} y_{-\beta_i}\right)$, with $\beta_1$,\ldots, $\beta_p$ the roots in $\Inv_{w}$ and
each $y_{\alpha}$ an element in the subspace of weight $\alpha$ 
(see \cite[Theorem 43]{BK} for the precise normalizations used in the definition of $\xi^{w}$).  
The factors of $(\frac{i}{2\pi})$ are taken care by the grading of the cohomology groups, and if ({\em v})
and ({\em ii}) hold we may also ignore the factors $\prod \langle\rho,\alpha\rangle$.
Thus an equivalent version of the above statements is 

\ResumeEnumerate
\item The map 
$$\phi'\colon \left(\H^{*}(\X,\QQ),\BKprod\right) \longrightarrow\left[\H^{*}(\gu^{+}_{\QQ})\otimes 
\H^{*}(\gu^{-}_{\QQ})\right]^{\gt}$$
defined by 
$$\phi'\left([\Omega_{w}]\right) = 
\left(-1\right)^{\frac{p(p-1)}{2}}
\xi^{w}$$
is also an isomorphism of graded rings, where, as above, $p=\ell(w)$.
\end{enumerate}

\np
Finally we note that the corresponding versions of these statements do not hold for the Belkale-Kumar
product on quotients $\G/\P$ in general.  For instance, when $\P$ is a maximal parabolic in type $\A$, 
the Belkale-Kumar product on $\H^{*}(\G/\P,\ZZ)$ is the usual cup-product, and there are many examples
of Littlewood-Richardson coefficients different from $0$ or $1$.

\bpoint{Acknowledgments} 
The method of \S\ref{sec:int-mult-one} using Weyl group combinatorics and representation theory is due to
P.\ Belkale and S.\ Kumar \cite{BK2}, and is used with their generous permission.
Ivan Dimitrov acknowledges excellent working conditions at the Max-Planck Institute.
Mike Roth acknowledges the hospitality of the University of Roma III.

\section{Intersection multiplicity one for classical groups}
\label{sec:int-mult-one}

\point
The main result of this note is the following theorem.

\tpoint{Theorem} \label{thm:mult1}
If $\G$ is classical (or $\G_2$) then condition \eqref{eqn:sym-liningup} implies that 
$\bigcap_{i=1}^{k} [\Omega_{w_i}]= 1$.

%
%
\np
We will compute the intersections by two different methods. In types $\A$, $\B$, and $\C$ we will use a
method combining Weyl group combinatorics and representation theory.  In type $\D$ we will use a more 
geometric fibration method, which however relies on a key combinatorial lemma.  
We now set up and apply the first method.

\bpoint{Torus fixed points and Weyl group combinatorics}\label{sec:torusfixed}
In this method we will compute the intersections by 
intersecting subvarieties representing these classes.   
The representatives will be torus stable subvarieties so it is useful to understand their torus fixed points.

\tpoint{Lemma} \label{lem:torusfixed}
For any element $w$ of $\W$ the torus fixed points of $(\wo w)^{-1}\Omega_{w} = (\wo w)^{-1} \X_{\wo w}$ 
are the elements of the set $$\left\{ u \st w \leq wu \rule{0cm}{0.6cm}\right\}.$$

\bpf
The torus fixed points of $\X_{\wo w}$ are the elements $v$ such that $v\leq \wo w$, and hence the torus fixed points
of $(\wo w)^{-1}\X_{\wo w}$ are the elements of the form $w^{-1}\wo v$ with $v\leq \wo w$.  

\np
Making the change of variables $u=w^{-1}\wo v$ (so that $v=\wo w u$), then this is the set of elements 
$\left\{u \st \wo w u \leq \wo w\right\}$. 
Since $\wo w u \leq \wo w$ if and only if $wu\geq w$ (in general $x\leq y$ iff $\wo x\geq \wo y$) 
this proves the lemma. \epf

\tpoint{Corollary} 
\label{cor:cycle-fixed-pts}
For any elements $w_1$, \ldots, $w_k\in \W$, the torus fixed points of the intersection 
$\bigcap_{i=1}^{k} (\wo w_i)^{-1} \Omega_{w_i}$ of the shifted Schubert varieties are the elements of the set 

\begin{equation}\label{eq:intersectiontorusfixed}
\left\{ u\in\W \bst w_i\leq w_iu\,\,\, \mbox{for all $i=1,\ldots, k$\rule{0cm}{0.6cm}}\right\}.
\end{equation}

\np
The proof of \cite[Lemma (2.6.1)]{DR} shows
that if $w_1$,\ldots, $w_k$ satisfy \eqref{eqn:sym-liningup} then
the intersection $\bigcap_{i=1}^{k} (\wo w_i)^{-1} \Omega_{w_i}$ is transverse at $e$, 
and that $e$ is an isolated component of the intersection.
Since the schemes $(\wo w)^{-1}\Omega_{w_i}$ are all fixed by the torus, any component of their intersection
must have a torus fixed point.  
Combining this with
Corollary \ref{cor:cycle-fixed-pts}, 
to prove Theorem \ref{thm:mult1} it is therefore sufficient (assuming \eqref{eqn:sym-liningup}) to show that

\begin{equation}\label{eqn:torus-fixed-e}
\left\{ u\in\W \bst w_i\leq w_iu\,\,\, \mbox{for all $i=1,\ldots, k$\rule{0cm}{0.6cm}}\right\} 
= \left\{{ e \rule{0cm}{0.5cm}}\right\}.
\end{equation}

\np
In order to demonstrate \eqref{eqn:torus-fixed-e} we will use the following well-known result (see for example
\cite[Theorem 7.7.7(i), p.\ 267]{dix}).

\tpoint{Proposition} \label{prop:order}
Let $x$, $y$ be elements of $\W$ with $x\leq y$ in the Bruhat order.  Then for any dominant weight $\lambda$ the
difference $x\lambda-y\lambda$ is a nonegative sum of positive roots.

\np
\tpoint{Lemma} \label{lem:combproof}
Suppose that $w_1$,\ldots, $w_k$ satisfy condition \eqref{eqn:sym-liningup}. 
Then 

\AlphaList
\begin{enumerate}
\item for each root $\alpha\in \Delta^{\!+}$ there is a  $w_i$ such that $w_i\alpha$ is a negative root.
\PauseEnumerate

\np
Further suppose that
$u$ is a solution to $w_i\leq w_i u$ for $i=1\ldots k$.  
Then for 
any dominant weight $\lambda$ we have :

\ResumeEnumerate
\item $\mu_{\lambda}:=\lambda-u\lambda$ is a nonnegative sum of positive roots. 
\item $w_i\mu_{\lambda}$ is a nonnegative sum of positive roots for $i=1$,\ldots, $k$. 
\item $\mu_{\lambda}$ is not a root or a multiple of a root.
\end{enumerate}

\bpf
Part ({\em a}) is obvious from condition \eqref{eqn:sym-liningup}. Part ({\em b}) follows from Proposition \ref{prop:order} 
since $u\geq e$ for any $u\in \W$. Part ({\em c}) follows from Proposition \ref{prop:order}, the condition that 
$w_i\leq w_i u$, and the obvious identity $ w_i\mu_{\lambda} = w_i\lambda-w_i u\lambda$.
Part ({\em d}) is proved by combining ({\em a}) and ({\em c}). \epf

\bpoint{Proof of Theorem \ref{thm:mult1} in types $\boldmath{\A}$, $\boldmath{\B}$, $\boldmath{\C}$}
Our strategy to show that \eqref{eqn:torus-fixed-e} holds (and thus that Theorem \ref{thm:mult1} holds),
is to assume that there is an element $u\neq e$ satisfying $w_i\leq w_iu$ for
$i=1,\ldots, k$ and then produce a dominant $\lambda$ such that $\mu_{\lambda}$ violates 
Lemma \ref{lem:combproof}({\em d}).  We now do this on a case-by-case basis.

\label{case:An}
\np
{\bf Type $\boldmath{\A_n}$.} 
Let $\ep_1$, \ldots, $\ep_{n+1}$ be a basis for the permutation
representation of $\W=\mathrm{S}_{n+1}$, where as usual the positive roots are of the form $\ep_{p}-\ep_{q}$ with $p<q$.
The fundamental weights  are
$\chi_p:=\ep_1+\ep_2+\cdots+\ep_p$ for $p=1$, \ldots, $n$. 
Let $u\in \W$ be such that $w_i\leq w_iu$ for $i=1$,\ldots, $k$.

\np
If $u\neq e$ then let $p$
be the smallest element of $\{1,\ldots, n\}$ such that $u\ep_p\neq\ep_p$.   
Since $p$ is the smallest such element, $u\ep_j=\ep_j$ for $j<p$ and $u\ep_p=\ep_q$ with $q>p$ and
hence $\mu:=\chi_p-u\chi_p=\ep_p-\ep_q$ is a positive root,
contradicting Lemma \ref{lem:combproof}({\em d}).  Therefore $u=e$ is the only possibility.

\label{case:Bn}
\np
{\bf Type $\boldmath{\B_n}$.} 
Let $\ep_1$, \ldots, $\ep_n$ be the usual basis upon which $\W$ 
operates by signed permutations.  The positive roots are $\ep_1$, \ldots, $\ep_n$ and elements of the form 
$\ep_p\pm\ep_q$ with $p<q$.
Fundamental
weights are $\chi_p=\ep_1+\cdots+\ep_p$ for $p=1,\ldots n-1$ and $\chi_n=\frac{1}{2}(\ep_1+\ep_2+\cdots+\ep_n)$.
Let $u\in \W$ be such that $w_i\leq w_iu$ for $i=1$,\ldots, $k$.

\np
If $u\neq e$ then let $p$ be the 
smallest element of $\{1,\ldots, n\}$ such that $u\ep_p\neq\ep_p$.   If $u\ep_p=\pm\ep_q$ with $q>p$ then 
$\mu:=\chi_p-u\chi_p=\ep_p\mp\ep_q$ is a positive root, contradicting Lemma \ref{lem:combproof}({\em d}).
Therefore $p=q$ and we must have 
$u\ep_p=-\ep_p$.  
But then $\mu:=\chi_p-u\chi_p$ is either twice a root (if $p<n$) or equal to a root (if $p=n$),
again contradicting Lemma \ref{lem:combproof}({\em d}).
Therefore $u=e$ is the only solution.

\label{case:Cn}
\np
{\bf Type $\boldmath{\C_n}$.} 
Let $\ep_1$, \ldots, $\ep_n$ be the usual basis upon which $\W$ 
operates by signed permutations.  The positive roots are $2\ep_1$, \ldots, $2\ep_n$ and elements of the form 
$\ep_p\pm\ep_q$ with $p<q$.  
Fundamental weights are $\chi_p=\ep_1+\cdots+\ep_p$ for $p=1,\ldots,  n$. 
Let $u\in \W$ be such that $w_i\leq w_iu$ for $i=1$,\ldots, $k$.

\np
The argument in this case is almost identical to that of $\B_n$: If $u\neq e$ let $p$ be the smallest element of
$\{1,\ldots, n\}$ such that $u\ep_p\neq\ep_p$.  Then $u\ep_p = \pm \ep_q$ with $q> p$ or $u\ep_p=-\ep_p$. In 
either case $\chi_p-u\chi_p$ is a root,
contradicting Lemma \ref{lem:combproof}({\em d}), so we again have $u=e$ as the only solution.

\np
Of course, the results for $\B_n$ and $\C_n$ are equivalent -- the natural isomorphism of Weyl groups respects the
Bruhat order, and induces a bijection of inversion sets (taking roots of the form $\ep_p\pm\ep_q$ to $\ep_p\pm \ep_q$
and roots of the form $\ep_p$ to $2\ep_p$); the proof above in the $\C_n$ case was included since it is equally short.

\np
We now turn to the general setup of the second method, which we will use in the proof of the theorem in type $\D$.

\bpoint{Restrictions of inversion sets and fibrations}\label{sec:restrictions-fibrations}
Recall that a subset $\SS\subseteq\delpos$ is called {\em closed} if whenever $\alpha,\beta\in\SS$ are such
that $\alpha+\beta$ is a root then $\alpha+\beta\in\SS$.  A subset $\SS$ is called {\em coclosed}\
if its complement $\SS^{\c}$ is closed.    We will use the following result of Kostant 
(see \cite[Proposition 5.10]{Ko}) without further reference: 
a subset $\SS$ of $\delpos$ is closed and coclosed if and only if $\SS=\invset_{w}$ for some $w\in\WW$;
the element $w$ is of course unique.

\tpoint{Definition}\label{def:phidef}
Let $\P\supseteq \B$ be a parabolic subgroup, $\Delta_{\P}$ the set of roots of $\P$, and $\W_{\P}$ 
the corresponding Weyl group, i.e., the group generated by the reflections in the roots contained of $\Delta_{\P}$, and
set $\Delta_{\P}^{+}=\Delta_{\P}\cap\Delta^{+}$.
For any $w\in \W$ the set $\invset_{w}\cap\delpos_{\P}$ is both closed and coclosed in $\delpos$, since 
$\invset_{w}$ is closed and coclosed in $\delpos$.  Therefore $\invset_{w}\cap\delpos_{\P}=\invset_{u}$ for a unique
$u\in \W$, and moreover $u\in \W_{\P}$.  We define $\phi_{\P}$ to be the (unique) map of sets
$\phi_{\P}\colon \W\longrightarrow \W_{\P}$ 
such that 

\begin{equation}\label{eqn:phidef}
\invset_{\phi_{\P}(w)} = \invset_{w} \cap \delpos_{\P}\quad\mbox{for all $w\in \W$.}
\end{equation}


\np
The map $\phi_{\P}$ has the following geometric meaning for projections of shifted Schubert varieties.

\tpoint{Proposition} \label{prop:parabolicproj}
Let $\P\supseteq \B$ be a parabolic subgroup, $\M=\G/\P$ and $\pi\colon\X\longrightarrow \M$ 
the projection.    Then for any $w\in \W$ :

\begin{enumerate}
\item $\pi(w^{-1}\X_{w})$ has dimension $\left|\invset_{w}\setminus\invset_{\phi(w)} \right|$.
\item Let $\G'$ be the Levi subgroup of $\P$ containing $\T$ and $\B':=\B\cap\G'$ the induced Borel.
As a subset of $\pi^{-1}(\pi(e))= \G'/\B'$, the fibre of $w^{-1}\X_{w}$ over $\pi(e)\in \M$ is 
$\phi(w)^{-1}\X_{\phi(w)}$. 
\end{enumerate}

\np
\bpf The composite $\B^{\op}\hookrightarrow \G\longrightarrow \G/\B$ is an open immersion of $\B^{\op}$ in $\X=\G/\B$. 
The image $\U$ ($\cong \B^{\op}$) of $\B^{\op}$ in $\X$ is therefore an affine space of dimension $\N$, 
whose torus-fixed coordinate rays are identified with the set $\Delta^{-}$ of negative roots.
Restricted to $\U$, each shifted Schubert variety $w^{-1}\X_{w}$ is the coordinate plane spanned by 
the coordinate vectors of the roots in $-\invset_{w}$.  The image of $\U$ in $\M$ 
is the affine space spanned by the roots in $\Delta^{-}\setminus \Delta_{\P}^{-}$, and the map $\pi$ restricted to 
$\U$ is the natural projection.  
The image of $\pi(w^{-1}\X_{w})$ restricted to $\pi(\U)$ is therefore the linear
space spanned by the roots in $-\invset_{w}\setminus\Delta_{\P}$, and the fibre in $\U$ over $\pi(e)$ is the 
linear space spanned by $-\invset_{w}\cap\Delta_{\P}=-\invset_{\phi(w)}$.  
This establishes both ({\em a}) and ({\em b}). \epf

\tpoint{Corollary} \label{cor:generic-X-proj-fibre}
For any $w\in \W$ the generic fibre of $\pi|_{\X_{w}}\colon\X_{w}\longrightarrow \pi(\X_{w})$ 
is $\X_{\phi(w)}$.

\np
\bpf
Since $\B$ acts transitively on an open subset of $\X_{w}$ containing $w\in \X$ 
it also acts transitively on an open subset of $\pi(\X_{w})$ containing $\pi(w)$.    Hence all fibres in this 
open set are isomorphic, and by Proposition \ref{prop:parabolicproj}({\em b}) the fibre over $\pi(w)$ is (after
shifting back) isomorphic to $\X_{\phi(w)}$. \epf

\np
We will also use the results above in the $\{[\Omega_{w}]\}_{w\in\W}$ basis:

\tpoint{Proposition} \label{prop:omegaproj}
With notation as above, for any $w\in \W$, 

\begin{enumerate}
\item $\pi(\Omega_{w})$ has codimension $\left|\invset_{w}\setminus\invset_{\phi(w)} \right|$ in $\M$.
\item The fibre of $(\wo w)^{-1}\Omega_{w}$ over $\pi(e)\in \M$ is $\phi(\wo w)^{-1}\Omega_{\phi(w)}$.
\item The general fibre of $\pi|_{g\Omega_{w_i}}\longrightarrow \pi(g\Omega_{w_i})$ 
is of the class $[\Omega_{\phi(w_i)}]$.
\end{enumerate}

\np
\bpf Parts ({\em a}) and ({\em b}) are restated versions of \ref{prop:parabolicproj}({\em a}) and ({\em b}), and
({\em c})  is a restated version of Corollary \ref{cor:generic-X-proj-fibre}. \epf

\np
Finally, we will use the following consequence of Proposition \ref{prop:omegaproj}({\em b}):

\tpoint{Corollary} \label{cor:inductive-intersect}
Suppose that $w_1$, \ldots, $w_k$ are elements of $\W$ such that $\sum \ell(w_i)=\N$, and 
let $\P$ be a parabolic subgroup, $\M=\G/\P$, and $\pi\colon \X\longrightarrow\M$ the projection.  We further
assume that $\sum_{i} |\invset_{\phi_{\P}(w_i)}| = |\Delta^{+}_{\P}|$ (note that 
by Proposition \ref{prop:omegaproj}({\em a})  
this is equivalent to the condition $\sum_i \codim(\pi(\Omega_{w_i}),\M)=\dim(\M)$).
Then we have the following equality of intersection numbers:

\begin{equation}\label{eqn:fibration-intersection}
\bigcap_{i=1}^{k} [\Omega_{w_i}] = 
\left({
\bigcap_{i=1}^{k} [\pi(\Omega_{w_i})]  
}\right) \cdot
\left({
\bigcap_{i=1}^{k} [\Omega_{\phi(w_i)}]  
}\right)
\end{equation}

\np
where the intersection on the left takes place in the cohomology ring $\H^{*}(\X,\ZZ)$, and the intersections on
the right in $\H^{*}(\M,\ZZ)$ and $\H^{*}(\X',\ZZ)$ respectively. 

\np
\bpf 
This is a consequence of Kleiman's transversality Theorem \cite[Corollary 4(ii)]{K}.
For general 
$g_1$, \ldots, $g_k$ in $\G$ the intersection of varieties $\bigcap_{i} g_i\Omega_{w_i}$ is transverse and computes
the intersection number $\bigcap_{i} [\Omega_{w_i}]$, and the points of intersection occur in the open cells of each
$g_i\Omega_{w_i}$.
We can also (by generality) choose the elements $g_i$ so that same holds for the intersection 
$\bigcap_{i=1}^{k} g_i\pi(\Omega_{w_i})$ in $\M$.

\np
The general fibre of $\pi|_{g\Omega_{w_i}}\longrightarrow \pi(g\Omega_{w_i})$ 
is of the class $[\Omega_{\phi(w_i)}]$ by Proposition 
\ref{prop:omegaproj}({\em c}).  Therefore, for
each point $p$ in the intersection $\bigcap_{i} g_i\Omega_{w_i}$, the fibre $\pi^{-1}(\pi(p))$ contains 
$\bigcap_{i=1}^{k} [\Omega_{\phi(w_i)}]$ points of the intersection.  Each of these
projects onto an intersection point of $\bigcap g_i\pi(\Omega_{w_i})$, which also meet transversely by our choice 
of $g_i$.
\epf

\bpoint{Description of fibration method} \label{sec:outline}
Assume that $w_1$,\ldots, $w_k$ satisfy \eqref{eqn:sym-liningup}.  
Let $\P \supset \B$ be a parabolic, $\M=\G/\P$, and $\pi\colon \X\longrightarrow \M$ the projection.
We also let $\G'$ be the Levi subgroup of $\P$, $\B'=\G'\cap \B$ the induced Borel, and $\X'=\G'/\B'$ the quotient.
Finally, let $\phi_{\P}\colon\W\longrightarrow\W_{\P}$ be the map of Definition \ref{def:phidef}.

\np
Since $\invset_{\phi(w_i)}=\invset_{w_i}\cap \Delta_{\P}$, if $w_1$, \ldots, $w_{k}$ satisfy 
\eqref{eqn:sym-liningup} then $\phi(w_i),\ldots, \phi(w_{k})$ satisfy 

\begin{equation}\label{eqn:parabolic_liningup}
\delpos_{\P}=\bigsqcup_{i=1}^{k} \invset_{\phi(w_i)}.
\end{equation}

\np
Thus we may apply Corollary \ref{cor:inductive-intersect}
to get the equality \eqref{eqn:fibration-intersection} of intersection numbers.
Suppose we can show that $\bigcap_{i=1}^{k} [\pi(\Omega_{w_i})] =1$ in $\H^{*}(\M,\ZZ)$, then 
\eqref{eqn:fibration-intersection} becomes
$ \bigcap_{i=1}^{k} [\Omega_{w_i}] = \bigcap_{i=1}^{k} [\Omega_{\phi(w_i)}]$. Since this second intersection is
taking place in $\H^{*}(\X',\ZZ)$, and since \eqref{eqn:parabolic_liningup} is simply condition 
\eqref{eqn:sym-liningup} for $\G'$, we may hope that we already know that the second intersection is $1$ by induction on
rank.  

\np
Thus the key inductive step for the fibration method is being able to show that the appropriate intersection 
in $\H^{*}(\M,\ZZ)$ is $1$.

\bpoint{Proof of Theorem \ref{thm:mult1} in type $\boldmath{\D}$ : Preliminaries} \label{sec:Dnstart} \label{sec:Dnroot}
We first prove, by induction, a combinatorial lemma 
(Lemma \ref{lem:Dn-combinatorial}).  This lemma and an elementary observation about the cohomology ring of
quadrics will establish the inductive step necessary to use the fibration method.

\np
{\bf $\D_n$ root systems.}
Let $\ep_1$, \ldots $\ep_n$ be the usual basis upon which $\W$ 
operates by signed permutations with an even number of sign changes.  The positive roots are elements of 
the form $\ep_p\pm\ep_q$ with $p<q$.
The fundamental
weights are $\chi_p=\ep_1+\cdots+\ep_p$ for $1\leq p\leq n-2$, 
$\chi_{n-1}=\frac{1}{2}(\ep_1+\ep_2+\cdots+\ep_{n-1}-\ep_n)$, and 
$\chi_{n}=\frac{1}{2}(\ep_1+\ep_2+\cdots+\ep_{n-1}+\ep_n)$.
We will also use $\D_2$ ($\cong \A_1\times \A_1$) and $\D_3$ ($\cong \A_3$) 
for the root systems defined as above with $n=2,3$.

\bpoint{Reductions to $\D_{n-1}$}
For any $p\in \{1,\ldots, n\}$, the subset of positive roots not involving $\ep_p$, i.e., the set 
$\delpos_p:=\{\ep_r\pm \ep_q\st r<q \,\mbox{and}\, r,q\neq p\}$, forms the positive roots of a 
sub-root system of type $\D_{n-1}$. 
For any element $w\in \W$, the intersection $\invset_{w}\cap \delpos_p$ is both closed and coclosed in $\delpos_p$, 
and hence is the inversion set of an element $\overline{w}$ in the $\D_{n-1}$ Weyl group.  
This reduction map from $\W_n$
to  $\W_{n-1}$ is not one coming from a parabolic $\P$ as in Definition \ref{def:phidef}, unless $p=1$.
Nonetheless, the reduction map exists and by definition has the property that if $w_1$,\ldots, $w_k$ satisfy 
\eqref{eqn:sym-liningup} then $\overline{w}_1$,\ldots, $\overline{w}_k$ also satisfy \eqref{eqn:sym-liningup} (i.e., 
$\delpos_p = \sqcup_{i=1}^{k} \invset_{\overline{w}_i}$).

\np
We will use this reduction map (the process of `deleting' an $\ep_p$) repeatedly in the proof of the combinatorial
lemma, so it is useful to understand the reduction explicitly.  We identify $\delpos_p$ with the $\D_{n-1}$ root
system on the basis elements $\overline{\ep}_1$,\ldots, $\overline{\ep}_{n-1}$ via the natural projection induced
by the linear map 
$$
\ep_q \longrightarrow 
\left\{
\begin{array}{cl}
\overline{\ep}_q & \mbox{if $q<p$ } \\
\overline{\ep}_{q-1} & \mbox{if $q>p$}. \\
\end{array}
\right.
$$

\np
Ignoring the signs for a moment, if we let $w$ act on $\ep_1$,\ldots, $\ep_{p-1}$, $\ep_{p+1}$, \ldots, $\ep_n$, then 
the order on the indices of resulting basis elements defines a permutation of $n-1$ objects.
The idea for the reduction $w\longrightarrow \overline{w}$ is that, treating $w$ and $\overline{w}$ as a signed 
permutations, the result of acting by $\overline{w}$ on 
$\overline{\ep}_1$, \ldots, $\overline{\ep}_{n-1}$
should induce the same relative order on the images as $w$ does above, and the signs should also be the same, with
the exception of the sign of $\overline{\ep}_{n-1}$, which may have to be switched to ensure an even number of total 
sign changes (i.e., if $w$ sends $\ep_p$ to the negative of some basis vector).  

\newcommand{\sgn}{\operatorname{sgn}}

\np
Explicitly, if $w(\ep_p) = \ep_{p'}$ for some $p'$ (as opposed to $w(\ep_p)=-\ep_{p'}$), then

$$
\overline{w}(\overline{\ep}_{q}) = \left\{
\begin{array}{ll}
\pm\overline{\ep}_{q'} & \mbox{if $q<p$, $w(\ep_q)=\pm\ep_{q'}$, and $q'<p'$} \\
\pm\overline{\ep}_{q'-1} & \mbox{if $q<p$, $w(\ep_q)=\pm\ep_{q'}$, and $q'\geq p'$} \\
\pm\overline{\ep}_{q'} & \mbox{if $q\geq p$, $w(\ep_{q+1})=\pm\ep_{q'}$, and $q'<p'$} \\
\pm\overline{\ep}_{q'-1} & \mbox{if $q\geq p$, $w(\ep_{q+1})=\pm\ep_{q'}$, and $q'\geq p'$} \\
\end{array}
\right.
$$

\np
where (for example in the first case above), the instructions mean 
$+\overline{\ep}_{q'}$ if $w(\ep_q)=+\ep_{q'}$ 
and
$-\overline{\ep}_{q'}$ if $w(\ep_q)=-\ep_{q'}$.

\np
If instead $w(\ep_p)=-\ep_{p'}$ for some $p'$ then $\overline{w}$ is the composite of the rule above 
followed by the map sending $\overline{\ep}_{n-1}$ to $-\overline{\ep}_{n-1}$ and acting as the identity on each 
$\overline{\ep}_i$, for $i< n-1$ (this ensures an even number of sign changes). 

\np
For our inductive argument we will need one fact which follows from the explicit formulae for the reduction. 
Suppose that $w$ is an element of $\W_{n}$ such that $w(\ep_1)=+\ep_q$ for $q<n$.  
If $\overline{w}$ is the result of deleting
some $\ep_p$ with $p>1$ then the only way that $\overline{w}(\overline{\ep}_1) = -\overline{\ep}_{q'}$ for some $q'$, 
or $\overline{w}(\overline{\ep}_1)=\overline{\ep}_{n-1}$ is if $w(\ep_1)=+\ep_{n-1}$ (i.e., $q=n-1$) 
and $w(\ep_p)=\pm\ep_n$. In this case $\overline{w}(\overline{\ep}_1)=\pm\overline{\ep}_{n-1}$, although we will not need this detail.

\np
We now prove the main combinatorial lemma for dealing with the $\D_n$ case.

\tpoint{Lemma} \label{lem:Dn-combinatorial}
Suppose that we are in the $\D_n$ case, and that $w_1$, \ldots, $w_k$ satisfy \eqref{eqn:sym-liningup}.  
Then there exists $i$ such that
$w_i(\ep_1)\in \{-\ep_1,-\ep_2,\ldots, -\ep_n, \ep_n\}$, i.e., for this $i$ either $w_i(\ep_1)=-\ep_p$ for some $p$, or
$w_i(\ep_1)=\ep_n$.

\bpf
Suppose that there is a counterexample for $\D_n$, i.e., $w_1$,\ldots, $w_k$ satisfying \eqref{eqn:sym-liningup} 
such that $w_i(\ep_1)=\ep_{p_i}$, $1\leq p_i \leq n-1$ for all $i=1$,\ldots, $k$.    We will show that we can always 
reduce
such a counterexample in $\D_{n}$ to a counterexample in $\D_{n-1}$.  For $n\geq 5$ this will follow by a counting
argument (and contradiction), for $n=3,4$ by a more detailed argument (and contradiction).  
Finally, it is obvious for $\D_2$ that no such counterexample exists, and this final contradiction proves the lemma.

\noindent
\underline{Case $\D_n$, $n\geq 5$.} 
We look for $\ep_p$, $2\leq p \leq n$ that we can `delete', and still maintain the 
counterexample.   
If it is not
possible to delete some $\ep_p$ and still maintain the counterexample, then for each $p$, $2\leq p\leq n$, 
there must be an $i_p$ such that $w_{i_p}(\ep_1)=\ep_{n-1}$ and 
$w_{i_p}(\ep_p)=\pm \ep_{n}$.  
The element $w_{i_p}$ then inverts exactly 
$n-2$ positive roots involving $\ep_1$ (exactly half of the positive roots $\{\ep_1\pm \ep_q \st q\neq p\}$).  
For different $p$, the corresponding $w_{i_p}$ are also distinct, since (for instance) $w_{i_p}^{-1}(\ep_{n})=
\pm \ep_{p}$.  Hence by 
\eqref{eqn:sym-liningup} 
these elements invert $(n-1)(n-2)$ distinct positive roots involving $\ep_1$.  Since there are exactly
$2(n-1)$ such roots, this gives the inequality $(n-1)(n-2)\leq 2(n-1)$ or $n\leq 4$.  Thus if $n\geq 5$ there is 
always such an $\ep_p$, and we can reduce the counterexample to the $\D_4$ case.

\np
\underline{Case $\D_4$.}  If there is no $p\in\{2,3,4\}$ so that we can delete $\ep_p$ and preserve the counterexample,
then as above there must be (reordering the $w_i$ as necessary) $w_1$, $w_2$, $w_3$ such that 
$w_i(\ep_1)=\ep_3$ and $w_i(\ep_{i+1})=\pm \ep_4$, $i=1,2,3$.  Each such $w_i$ inverts exactly two roots involving
$\ep_1$, and hence we must have $w_i(\ep_1)=\ep_1$ for all $i\geq 4$ (if $k\geq 4$), since there are exactly
six positive roots of the form $\ep_1\pm \ep_q$, $q\in \{2,3,4\}$.

\np
From the conditions, $w_1$ inverts exactly one of $\ep_2\pm\ep_3$ and exactly one of $\ep_2\pm \ep_4$.  This implies
that $w_2(\ep_2)=+\ep_q$ with $q\in \{1,2\}$, since if $w_2(\ep_2)=-\ep_q$ then both of $\ep_2\pm \ep_3$ would be
inverted by $w_2$, contradicting the fact that $w_1$ inverts exactly one of them, and condition \eqref{eqn:sym-liningup}.
Similarly, we must have $w_3(\ep_2)=+\ep_q$ with $q\in\{1,2\}$ or $w_3$ would invert both of $\ep_2\pm\ep_4$.

\np
But now none of $w_1$, $w_2$, and $w_3$ inverts $\ep_1+\ep_2$, and since $w_i(\ep_1)=\ep_1$ for all $i\geq 4$, we
see that $\ep_1+\ep_2$ is never inverted, again contradicting \eqref{eqn:sym-liningup}.  Thus we may reduce the 
counterexample to the $\D_3$ case.

\np
\underline{Case $\D_3$.} Again, assume that there is no $p\in\{2,3\}$ which can be deleted and maintain the counterexample.
Then (after reordering) we must have $w_i(\ep_1)=\ep_2$, $w_i(\ep_{i+1})=\pm\ep_3$ for $i=1,2$, and $w_i(\ep_1)=\ep_1$
for $i\geq 3$.  Again $w_1$ inverts exactly one of $\ep_2\pm\ep_3$, so condition \eqref{eqn:sym-liningup} implies that
we must have $w_2(\ep_2)=\ep_1$.   But now, as before, no $w_i$ inverts $\ep_1+\ep_2$, a contradiction.  Thus
we can reduce any counterexample in $\D_3$ to $\D_2$.

\np
\underline{Case $\D_2$.} The condition for the counterexample now means that $w_i(\ep_1)=\ep_1$ for all $i$, and hence
$w_i(\ep_2)=\ep_2$ for all $i$ (since each $w_i$ is a signed permutation of $\ep_1$, $\ep_2$ with an even number
of sign changes).  I.e., each $w_i=e$.  This certainly contradicts \eqref{eqn:sym-liningup}, and hence no such
counterexample exists. 

\np
This finishes the proof of Lemma \ref{lem:Dn-combinatorial}. \epf

\bpoint{Intersections on even dimensional quadrics}\label{sec:intersection}
Let $\P_1$ be the parabolic so that $\W_{\P_1}\subset\W$ is the stabilizer of $\epsilon_1$.  Then $\Q_{n}:=\G/\P_1$ is a 
smooth quadric hypersurface in $\PP^{2n-1}$.  Let $\pi\colon \X \longrightarrow \Q_{n}$ be the projection.  
The fibre $\X'=\pi^{-1}(\pi(e))$ is of type $\D_{n-1}$. 

\np
The cohomology ring of $\Q_n$ is generated by $h$ (the class of a hyperplane section) and 
two classes $a$ and $b$ of codimension $(n-1)$ (i.e., in the middle cohomology), satisfying the relations 

\begin{equation}\label{eqn:Qnrelations}
\rule{1.1cm}{0cm}
\mbox{\raisebox{0.0cm}{$
h^{n-1}=a+b,ha=hb, h^na=0, a^2=b^2=\frac{1}{2}(1-(-1)^{n})[pt], ab=\frac{1}{2}(1+(-1)^{n})[pt]
$}}
\end{equation}

\np
where $[pt]$ is the class of a point.  The classes $a$ and $b$ are represented by linear subspaces of $\PP^{2n-1}$
of dimension $n-1$ contained in $\Q_{n}$.

\np
The cohomology ring therefore has the presentation

$$
\H^{*}(\Q_{n},\ZZ) = \frac{\ZZ[h,a,b]}{(\mbox{relations in \eqref{eqn:Qnrelations}})}.
$$

\np
We will use the integral basis for $\H^{*}(\Q_n,\ZZ)$ given by $\{h^{k}\}_{0\leq k\leq n-2}$ in codimension $\leq n-2$, 
$a$ and $b$ in codimension $n-1$, and $\{h^{k}a\}_{1\leq k\leq n-1}$ in codimensions $n$ to $2(n-1)$.
Under the projection $\pi$, the image of each Schubert cell in $\X$ is sent to a variety whose cohomology class
is one of the integral basis classes above.  The complex codimension of the image of $\Omega_{w}$ is 
the number of roots involving $\ep_1$ (the roots of the form $\ep_1\pm\ep_q$) in $\invset_{w}$.

\np
Since $\Q_{n}$ has degree $2$, $h^{2n-2}=2[pt]$, and since $a$ is the class of a linear space $h^{n-1}a=1[pt]$.
Given our choice of basis classes, this immediately proves the following result.

\tpoint{Lemma} \label{lem:intersection}
Let $c_1$,\ldots, $c_k$ be basis cohomology classes in $\H^{*}(\Q_{n},\ZZ)$ whose (complex) 
codimensions sum to $2(n-1)=\dim(\Q_n)$.  Then 

$$\bigcap_{i=1}^{k} c_i = \left\{
\begin{array}{cl}
1 & \mbox{if some $c_i$ has codimension $\geq n-1$} \\
2 & \mbox{if all $c_i$ have codimension $\leq n-2.$} \\
\end{array}
\right.
$$

\bpoint{Proof of Theorem \ref{thm:mult1} in type $\boldmath{\D}$ : Geometric Approach} 
\label{sec:Dnstart} \label{sec:Dnroot}
\label{case:Dn}
We prove the result by induction on $n$.  The case $n=3$
is $\D_3=\A_3$, which is covered by \S\ref{case:An}.  It therefore suffices to give the inductive step.  Suppose
that $w_1$,\ldots, $w_k$ satisfy \eqref{eqn:sym-liningup}, then by Lemma \ref{lem:Dn-combinatorial} there is some
$i$ so that $w_i(\ep_1)\in\{-\ep_1,\ldots, -\ep_n,\ep_n\}$.  For such an $i$, $w_i$ inverts at least $n-1$ positive
roots involving $\ep_1$ (the roots of the form $\ep_1\pm\ep_q$).  Hence $\pi(\Omega_{w_i})$ has codimension $\geq n-1$
in $\Q_{n}$ by Proposition \ref {prop:omegaproj}({\em a}).
By Lemma \ref{lem:intersection} this means that
$\bigcap_{i=1}^{k} [\pi(\Omega_{w_i})]=1$, and hence by Corollary \ref{cor:inductive-intersect} that 

$$ 
\bigcap_{i=1}^{k} [\Omega_{w_i}] = 
\bigcap_{i=1}^{k} [\Omega_{\phi(w_i)}],
$$

\np
where 
$\phi$ is the map $\phi\colon\W\longrightarrow\W_{\P_1}$ of Definition \ref{def:phidef}.  Since $\phi(w_1)$,
\ldots, $\phi(w_k)$ are elements of the $\D_{n-1}$ root system satisfying \eqref{eqn:sym-liningup}, we conclude by
the inductive hypothesis that we have $\bigcap_{i=1}^{k} [\Omega_{w_i}] = 1$.  \epf

\np
It is also possible to use the method of \S\ref{sec:torusfixed} to prove Theorem \ref{thm:mult1} in the $\D_n$ case;
the key step is again Lemma \ref{lem:Dn-combinatorial}.  To avoid some extra combinatorial digressions, we only
sketch the argument.

\bpoint{Proof of Theorem \ref{thm:mult1} in type $\boldmath{\D}$ : Combinatorial Approach} 
Let $u$ be such that $w_i\leq w_i u$ for $i=1$,\ldots, $k$.  We want to 
show that $u=e$.   We first show that $u(\ep_1)=\ep_1$.  If $u(\ep_1)=\pm\ep_q$ with $q>1$ then $\mu=\chi_1-u\chi_1=
\ep_1\mp\ep_q$ is a positive root, contradicting Lemma \ref{lem:combproof}({\em d}).  If $u(\ep_1)=-\ep_1$, then
$\mu:=\chi_1-u\chi_1=2\ep_1$, which is a sum of positive roots.  However, by Lemma \ref{lem:Dn-combinatorial}
there is an $i$ such that $w_i(\ep_1)=-\ep_p$ or $w_i(\ep_1)=\ep_n$.  Then $w_i\mu=-2\ep_p$ or $w_i\mu=2\ep_n$, 
neither of which are sums of positive roots, contradicting Lemma \ref{lem:combproof}({\em c}).
Thus we must have $u(\ep_1)=\ep_1$, and so $u\in \W_{\P_1}$.  Applying the map $\phi\colon\W\longrightarrow\W_{\P_1}$
one can check (these details are omitted) that $\phi(w_i)\leq \phi(w_i)\phi(u)$, where the order is
now the Bruhat order on $\W_{\P_1}$.   By induction, the only solution is $\phi(u)=e$, and since $u\in\W_{\P}$,
this implies that $u=e$.  \epf

\bpoint{Proof of Theorem \ref{thm:mult1} in type $\boldmath{\G_2}$} 
The argument is elementary for any rank $2$ group; by Lemma
\ref{lem:maxk} below there are at most two $w_i$ with $w_i\neq e$.  If there are exactly two such $w_i$, say $w_1$ 
and $w_2$ then condition \eqref{eqn:sym-liningup} implies that $\Omega_{w_1}$ and $\Omega_{w_2}$ are Poincar\'e dual
pairs, so 
$\bigcap_{i=1}^{k} [\Omega_{w_i}] = [\Omega_{w_1}]\cap[\Omega_{w_2}] =1$.  If there is only one such $w_i$,
then it must be equal to $\wo$ and since $\Omega_{\wo}$ is a point we again have $\bigcap_{i=1}^{k} [\Omega_{w_i}] = 1$.

\tpoint{Lemma}\label{lem:maxk} If $w_1$, $\ldots$, $w_{k}\in \W$ satisfy \eqref{eqn:sym-liningup}
then the number of $w_i$ which are not equal 
to $e$ is at most $\rank(\G)$.

\np
\bpf
Each $\invset_{w_i}$ is coclosed, so if $\invset_{w_i}$ does not contain any simple roots, then $\invset_{w_i}$ 
does not contain any roots and therefore $\invset_{w_i}=\emptyset$ and so $w_i=e$.    Therefore if $w_i\neq e$
the set $\invset_{w_i}$ contains a simple root.  Since the union $\delpos=\sqcup_{i} \invset_{w_i}$ is disjoint,
the number of $w_i$ with $w_i\neq e$ is therefore at most the number of simple roots. 
\epf

\section{Appendix : The equivalence of conditions ({\em iii}) and ({\em iv})}
\label{sec:Appendix}

\noindent
It is clear that ({\em iii}) (being the case $k=3$ of ({\em iv})) is implied by ({\em iv}).
To prove the other direction we will need to discuss the product $\BKprod$ in more detail.
This product is obtained by specializing a deformation of the ordinary cup product.
This deformation was introduced by Belkale and Kumar.

\bpoint{The Belkale-Kumar deformation of the cup product on $\G/\B$} Let $\alpha_1$,\ldots, $\alpha_{n}$ denote
the simple roots of $\G$ and let $\sQ$ be the root lattice.
Introduce variables $\tau_1$,$\tau_2$,\ldots, $\tau_n$, one for each simple root.  
For any $\gamma\in \sQ$ we use the notation $\taubar^{\gamma}$ to denote the Laurent monomial 
$\tau_1^{m_1}\tau_2^{m_2}\cdots\tau_n^{m_n}$,
where $\sum_{i=1}^{n} m_i\alpha_i = \gamma$ is the unique expression of $\gamma$ as a $\ZZ$-linear combination 
of simple roots.

\np
Following \cite[Definition 5]{BK} for any $w\in \W$ we define 
$\chi_{w}=\sum_{\alpha\in \Inv_{w}} \alpha$.
The operation $\odot$ acting on two basis classes is defined 
\cite[p.\ 199]{BK} by the formula 
$$
[\Omega_{w_1}]\odot [\Omega_{w_2}] : = 
\sum_{w} \taubar^{(\chi_{w}-\chi_{w_1}-\chi_{w_2})} c_{w_1,w_2}^{w} [\Omega_{w}].
$$
Belkale and Kumar \cite[Proposition 17(a)]{BK} prove that if $c_{w_1,w_2}^{w}\neq 0$ then 
$\chi_{w}-\chi_{w_1}-\chi_{w_2}$ is in the positive root lattice, and thus all exponents of
$\taubar^{(\chi_{w}-\chi_{w_1}-\chi_{w_2})}$ are nonnegative. The product above therefore takes 
values in $\H^{*}(\X,\ZZ)\otimes \ZZ[\tau_1,\ldots, \tau_n]$.  The operation $\odot$ is then 
extended to all of $\H^{*}(\X,\ZZ)\otimes \ZZ[\tau_1,\ldots,\tau_n]$ by $\ZZ[\tau_1,\ldots,\tau_n]$-linearity.  

\np
From the formula it is clear that $\odot$ is commutative.
One checks by induction (see \cite[Proposition 17(c)]{BK}) that for any $w_1,\ldots, w_k\in \W$ 

\begin{equation}\label{eqn:general-BK-deform}
[\Omega_{w_1}]\odot [\Omega_{w_2}]\odot \cdots \odot [\Omega_{w_k}] = 
\sum_{w\in \W} \taubar^{(\chi_{w}-\sum\chi_{w_i})} c^{w}_{w_1,\ldots, w_k} [\Omega_{w}],
\end{equation}
where $c^{w}_{w_1,\ldots, w_k}$ is the coefficient of $[\Omega_{w}]$ in the expression of 
$\cap_{i=1}^{k} [\Omega_{w_i}]$ as a sum of basis classes.
Thus the product of basis elements in the deformed product is the usual cup product, with each term in the result
shifted by a mononomial in $\tau_1$,\ldots, $\tau_n$, where the monomial depends on the term and the classes being 
multiplied.

\np
Setting all $\tau_i=1$ recovers the usual cup product. The Belkale-Kumar product $\BKprod$ is defined
as the specialization obtained by setting all $\tau_i=0$.

\tpoint{Lemma} \label{lem:unions-of-inv-sets} 
If $w_1$,\ldots, $w_k\in \W$ satisfy \eqref{eqn:lining-up}, then 

\AlphaList
\begin{enumerate}
\renewcommand{\itemsep}{0.25cm}
\item 
$[\Omega_{w_1}]\BKprod [\Omega_{w_2}]\BKprod \cdots \BKprod [\Omega_{w_k}] = 
[\Omega_{w_1}]\cap [\Omega_{w_2}]\cap \cdots \cap [\Omega_{w_k}].$

\item
For any subset $\I\subseteq\{1,2,\ldots, k\}$ there is an element $w\in \W$ so that
$\Inv_{w}=\sqcup_{i\in \I}\Inv_{w_i}$.   
\end{enumerate}

\bpf
One of the properties of the inversion sets is that 
for any $w\in \W$, $\ell(w)=|\Inv_{w}|$.  Therefore if $w_1$,\ldots, $w_k$ satisfy
\eqref{eqn:lining-up} we have $\sum_{i=1}^{k} \ell(w_i) = \sum_{i=1}^{k} |\Inv_{w_i}|=|\Delta^{+}|=\dim(\X)$.
The only class in dimension zero is the class of a point, $[\Omega_{\wo}]$.
Since (again by \eqref{eqn:lining-up}) we have $\chi_{\wo}=\sum_{i=1}^{k} \chi_{w_i}$, 
we conclude by \eqref{eqn:general-BK-deform} that 
$$
[\Omega_{w_1}]\odot [\Omega_{w_2}]\odot \cdots \odot [\Omega_{w_k}] = 
c^{\wo}_{w_1,\ldots, w_k} [\Omega_{\wo}] = \bigcap_{i=1}^{k} [\Omega_{w_k}].
$$
This proves ({\em a}).

\np
For a proof of part ({\em b}), see 
\cite[Corollary 5.4.9]{DR} or \cite[Proposition 2.7]{ddmrww} 
(the proof in the second reference is presented in a more combinatorial context, and in the language of type $\A$,
but works in all types).
Part ({\em b}) may also be deduced using Lie algebra cohomology.
\epf

\bpoint{Proof that ({iii}) implies ({iv})} By Lemma \ref{lem:unions-of-inv-sets}({\em a}) it is sufficent
to show that $[\Omega_{w_1}]\BKprod \cdots \BKprod [\Omega_{w_k}]=1$. 
By part ({\em b}) of the same lemma there is an element $u\in \W$ such that
$\Inv_{u}=\Inv_{w_{k-1}}\sqcup \Inv_{w_k}$.   By \eqref{eq:structure-coeffs} and ({\em iii}) (in its equivalent
form ({\em ii})) we have $[\Omega_{w_{k-1}}]\BKprod [\Omega_{w_k}]=1[\Omega_{u}]$.  Thus
$$
[\Omega_{w_1}]\BKprod[\Omega_{w_2}]\BKprod\cdots\BKprod 
[\Omega_{w_{k-2}}]\BKprod [\Omega_{w_{k-1}}] \BKprod [\Omega_{w_k}]=
[\Omega_{w_1}]\BKprod[\Omega_{w_2}]\BKprod  \cdots \BKprod [\Omega_{w_{k-2}}]\BKprod [\Omega_{u}]$$
with $(\bigsqcup_{i=1}^{k-2} \Inv_{w_i})\sqcup \Inv_{u}=\Delta^{+}$.
I.e., we have reduced the expression we are interested in to a similar expression with one fewer term.
Continuing in this manner we reduce the expression to $[\Omega_{\wo}]$, the class of a point. \epf



\begin{thebibliography}{[Bott]}

\bibitem[BK]{BK} P.\ Belkale and S.\ Kumar, {\em Eigenvalue problem and a new product in cohomology of flag varieties},
Invent. Math. {\bf 166} (2006), 185--228.

\bibitem[BK2]{BK2} P.\ Belkale and S.\ Kumar, {\em private communication}.









\bibitem[Dix]{dix} J.\ Dixmier, {\em Enveloping algebras}.  Revised reprint of the 1997 translation.  Graduate
Studies in Mathematics {\bf 11}.  American Mathematical Society, Providence, RI, 1996.

\bibitem[D-W]{ddmrww} R.\ Dewji, I.\ Dimitrov, A.\ McCabe, M.\ Roth, D.\ Wehlau, and J.\ Wilson,
{\em Decomposing inversion sets of permutations and applications to faces of the Littlewood-Richardson cone},
J. Algebraic Combin.\ {\bf 45}:4 (2017), 1173--1216.

\bibitem[DR]{DR} I.\ Dimitrov and M.\ Roth, {\em Cup products of line bundles on homogeneous varieties and
generalized PRV components of multiplicity one}, to appear in Algebra \& Number Theory.







\bibitem[K]{K} S.\ Kleiman, {\em The transversality of a general translate}, Compositio Math.\  {\bf 28} (1974),
 287--297. 


\bibitem[Ko]{Ko} B.\ Kostant, {\em Lie algebra cohomology and the generalized Borel--Weil theorem}, 
Ann. of Math. (2) {\bf 74} (1961), 329--387.
















\end{thebibliography}
\end{document}